\documentclass[12pt,a4paper,reqno]{amsart}

\usepackage{amsfonts}
\usepackage{graphicx}
\usepackage{graphics}
\usepackage{epsfig}
\usepackage{anysize}
\usepackage{lipsum}
\usepackage{wrapfig}
\marginsize{3cm}{3cm}{3cm}{3cm}
\newtheorem{theorem}{Theorem}[section]
\newtheorem{proposition}{Proposition}[section]
\newtheorem{lemma}{Lemma}[section]
\newtheorem{remark}{Remark}[section]
\newtheorem{corollary}{Corollary}[section]
\newtheorem{definition}{Definition}[section]
\newtheorem{example}{Example}[section]

\makeatletter

\DeclareMathOperator{\Image}{Im}
\DeclareMathOperator{\IMA}{im}
\DeclareMathOperator{\Rmo}{Rmo}
\DeclareMathOperator{\Plat}{Plat}
\DeclareMathOperator{\Fxdp}{Fix}
\DeclareMathOperator{\Exc}{Exc}
\DeclareMathOperator{\Wex}{Wex}
\DeclareMathOperator{\AX}{Ax}
\DeclareMathOperator{\AXL}{AxL}
\DeclareMathOperator{\exc}{exc}
\DeclareMathOperator{\fxdp}{fix}
\DeclareMathOperator{\plat}{plat}
\DeclareMathOperator{\wex}{wex}
\DeclareMathOperator{\ax}{ax}
\DeclareMathOperator{\nom}{nom}
\DeclareMathOperator{\F}{F}

\DeclareMathOperator{\flip}{Flip}

\begin{document}

	\setcounter{page}{1}
	\vspace{2cm}
	\author[]{Fufa Beyene$^1$, Roberto Mantaci$^2$ }
	\title[\centerline{Permutations with non-decreasing transposition array ...
		\hspace{0.5cm}}]{Permutations with non-decreasing transposition array and pattern avoidance}

	\thanks{\noindent $^1$ Addis Ababa University- Department of Mathematics - 1176 Addis Ababa - Ethiopia\\
		\indent \,\,\, e-mail: fufa.beyene@a{}au.edu.et\\ 
		\indent $^2$ Universit\'e de Paris- IRIF - 8 Pl. Aur\'elie Numours, 75013 Paris - France\\
		\indent \,\,\, e-mail: mantaci@irif.fr} 
	
	%
	\begin{abstract}
		We give some results about a bijection associating each permutation with a subexcedant function (i.\ e., a function $f : [n] \rightarrow [n]$ such that $f(i) \leq i, 1\le i\le n$). This function is related to a particular decomposition of the permutation as a product of transpositions and therefore it has been called transposition array in the literature. In particular, we identify anti-exceedance positions of the permutation through its transposition array and we give an expression of the bijection in terms of the cycle structure of the permutation. We give a characterization of a family of permutations having a non-decreasing transposition array and study length $3$ pattern avoidance therein.
		
		
		
	\end{abstract}
	\maketitle 
	\bigskip
	\section{Introduction}
	Let $n$ be a fixed positive integer and $[n]:=\{1,2,\ldots,n\}$. At least since the works of Lehmer (\cite{Le}, 1960) or perhaps even much earlier (Laisant \cite{La}, 1888), permutations over $[n]$ have been coded with integer $n$-tuples $(a_1,a_2,\ldots,a_n)$ such that $1 \leq a_i \leq i, \forall i \in [n]$, or equivalently with functions $f:[n]\longrightarrow[n]$ such that $1 \leq f(i) \leq i, \forall i \in [n]$, called \emph{subexcedant functions}. A subexcedant function $f$ is written as the word $f_1f_2\cdots f_n$, where $f_i=f(i), \forall i\in[n]$. 
	
	Permutation codes are interesting because certain algorithms perform better over the codes than they do over the permutations themselves. Codes allow for instance to implement efficient algorithms for the exhaustive generation of some specific classes of  permutations. To do so, one often has to ``read'' the properties of the permutation in its code and this gives birth to interesting combinatorial problems (see for instance \cite{Ba-Va,Bo,Fo-Ze,Ma-Ra,St1,Te}).
	
	Any permutation over $[n]$ can be decomposed in several ways as a product of \emph{transpositions}, cycles of length $2$, but it can be decomposed only in one way as a product of transpositions of the kind $(n,f_n)(n{-}1,f_{n-1})\cdots(1,f_1)$ and in one way as $(1,f_1')(2 ,f_{2}')\cdots(n,f_n')$, where $f_1f_2\cdots f_n=f$ and $f_1'f_2'\cdots f_n'=f'$ are subexcedant functions over $[n]$. In the product of these transpositions, if $f(i)=i$ for some $i$ correspondingly for $f'$, then the corresponding transposition degenerates into the identity.
	The first version, associating the code $f$ to the permutation $\sigma =\phi(f)= (n,f_n)(n{-}1,f_{n-1})\cdots(1,f_1)$ was studied for instance by Rakotondrajao and the second author (see \cite{Ma-Ra}), who related excedances to appropriate statistics on the corresponding codes. The second version, associating the code $f'$ to the permutation $\sigma =\tilde{\phi}(f)= (1,f_1')(2,f_2')\cdots(n,f_n')$ was independently introduced by Baril (see \cite{Ba1}), and who also studied in \cite{Ba} in particular the positions of weak excedances in a permutation using the corresponding code and who designates this permutation code with the term \emph{transposition array}, a name that can be sensibly attributed to the bijection of the first version as well.  
	
	In this article we study transposition arrays obtained by the bijection $\phi$ described in \cite{Ma-Ra}. The organization of the paper is as follows: In section 2, we recall some definition and some basic results about this permutation code. In section 3, we provide a new definition of the bijection $\phi$ based on the decomposition of the permutation as a product of disjoint cycles. In section 4, we study  permutations whose transposition array is non-decreasing. We conclude the article with a section studying the avoidance of length $3$ patterns in the set of permutations with non-decreasing transposition arrays.
	\section{Notation and preliminaries}
	Let
	$\mathfrak{S}_n$ be the symmetric group of permutations over the set $[n]$. A permutation $\sigma\in\mathfrak{S}_n$ can be given as a word $\sigma=\sigma(1)\sigma(2)\cdots\sigma(n)$
	or in cycle notation, where as usual a \emph{cycle} in $\sigma$ can be written as $(j,\sigma(j),\sigma^2(j), \ldots, \sigma^{t-1}(j))$, where $t$, the \emph{length} of the cycle, is the first positive integer such that $\sigma^t(j)=j$. Cycles of length one are \emph{fixed points}. The cycle notation is $\sigma=C_1C_2\cdots C_k$, where the $C_i$'s are the different and thus disjoint cycles of $\sigma$.
	
	A permutation $\sigma$ is said to have an excedance at $i\in[n]$ if $\sigma(i)> i$, a weak excedance if $\sigma(i)\ge i$, and an anti-excedance if $\sigma(i)\leq i$. If $\sigma$  has an excedance (respectively, weak excedance or anti-excedance) at $i$, $\sigma(i)$ is called an excedance (respectively, weak excedance or anti-excedance) letter. If $\sigma=\sigma(1)\sigma(2)\cdots\sigma(n)\in\mathfrak{S}_n$, then we use the notation: \begin{align*}
		&\Exc(\sigma):=\{i\in[n]: \sigma(i)> i\} \text{ and}\\
		&\Wex(\sigma):=\{i\in[n]: \sigma(i)\ge i\},\\
		&\AX(\sigma):=\{i\in[n]: \sigma(i)\leq i\}.
	\end{align*}
	Furthermore, $\exc(\sigma):=|\Exc(\sigma)|, \wex(\sigma):=|\Wex(\sigma)|$, and $\ax(\sigma):=|\AX(\sigma)|$.
	
	Recall that a subexcedant function is a function $f:[n] \rightarrow [n]$ such that $1\leq f_i\leq i, \forall i\in [n]$. We let $\F_n$ denote the set of all subexcedant functions over $[n]$. Let $f=f_1f_2\cdots f_n\in \F_n$. Then the image set $\Image(f)$ denote the set of elements in $f([n])$ and $\IMA(f)$ is its cardinality. 
	Further we define the \textit{statistics} $\Rmo$, $\Plat$ and $\Fxdp$ as follows:
	\begin{align*}&\Rmo(f):=\{i : 1\le i\le n, f_i\notin\{f_{i+1}, \ldots, f_n\}\},\\
		&\Plat(f):=\{i : 1\leq i< n, f_i=f_{i+1}\}, ~\plat(f):=|\Plat(f)|,\\
		&\Fxdp(f):=\{i : 1\leq i\leq n, f_i=i\}, ~\fxdp(f):=|\Fxdp(f)|.
	\end{align*}
	
	
	\subsection{The inverse function of $\phi$}
	As shown in \cite{Ma-Ra}, the procedure to obtain the subexcedant function $f=f_1f_2\cdots f_n$ associated with a permutation
	$\sigma=\sigma(1)\sigma(2)\cdots\sigma(n)$ such that  $\phi(f)=\sigma$ is described by an $n$-step loop.
	\begin{itemize}
		\item[--] For $i = n, n-1, \ldots, 1$: \vspace{-.1cm}
		\begin{itemize}	
			\item set $f_i=\sigma(i)$
			\item set $\sigma =(i, \sigma(i))\circ \sigma  $ (this swaps $i$ and $\sigma(i)$ in the permutation)
		\end{itemize}
	\end{itemize}
	For $i=n, n-1, \ldots, 1$ we let $\sigma_i$ be the permutation at the  iteration of the ``for'' loop corresponding to a given $i$,  then the procedure sets $f_{i}=\sigma_i(i)$, where $\sigma_n=\sigma$ and $\sigma_i= (i+1, \sigma_{i+1}(i+1))\circ\sigma_{i+1}$.
	
	\begin{example}
		Take $\sigma=23154$. Then $\sigma_5=2315\underline{4}$ and $f_5=\sigma_5(5)=\sigma(5)=4$ and $\sigma_4=(5,4)\circ\sigma_5=231\underline{4}\mathbf{5}$. Now $f_4=\sigma_4(4)=4$ and $\sigma_3=(4,4)\circ\sigma_4=23\underline{1}\mathbf{45}$. From $\sigma_3$ we have $f_3=\sigma_3(3)=1$ and $\sigma_2=(3,1)\circ\sigma_3=2\underline{1}\mathbf{345}$. So $f_2=\sigma_2(2)=1$ and $\sigma_1=(2,1)\circ\sigma_2=\underline{1}\mathbf{2345}$. Finally from $\sigma_1$ we have $f_1=1$. Therefore, $f=f_1f_2f_3f_4f_5=11144$.
	\end{example}
	\begin{remark}
		We note that at the end of the iteration for a given $i$, the permutation $\sigma_i$ fixes all points greater than $i$ (whence it may be regarded as a permutation on $[i]$). Therefore, the algorithm sorts the permutation $\sigma$ (by selecting the maximum and swapping).
	\end{remark}
	\subsection{Some properties of $\phi$}
	\subsubsection{The insertion method}\label{insertionmthd}
	By observing the sorting process of the permutation $\sigma$ when we compute its corresponding transposition array $f=\phi^{-1}(\sigma)$, we have another way to express and compute $\sigma=\phi(f)$ other than as a product of transpositions. Indeed, a given subexcedant function $f= f_1f_2\cdots f_n$ contains the information to easily construct  $\sigma$ by ``undoing'' that sorting process. The algorithm follows.
	Start with the identity $\sigma_1= 1$, suppose that the integers $1, 2, \ldots, i-1$ have already been inserted. Then replace the integer $f_i$ by $i$ and append the integer $f_i$ at the end. If $f_i=i$, then simply append $f_i$ at the end.
	\begin{example}
		\label{eg3.1}
		Take $f=11144\in \F_5$, then we construct  $\sigma=\phi(f)$ as follows.
		\begin{align*}
			&\sigma_1=1\\
			&\sigma_2=21\\
			&\sigma_3=231\\
			&\sigma_4=2314\\
			&\sigma_5=23154=\sigma
		\end{align*}
	\end{example}
	As shown in \cite[Proposition 3.5]{Ma-Ra}, the set $\Image(f)$ coincides with the set of anti-excedance letters of $\sigma=\phi(f)$, i.\ e., $i \in \Image(f)$ if and only if $\sigma^{-1}(i)\in\AX(\sigma)$. Independently, in \cite[Lemma 4]{Ba}, Baril has proved that the set $\Image(f)$ coincides with the set of weak excedances of $\sigma^{-1}=\tilde{\phi}(f)$. Actually these two statements are equivalent, since obviously we have
	\begin{lemma}
		The integer $i\in[n]$ is an anti-excedance of $\sigma$ if and only if $\sigma(i)$ is a weak excedance of $\sigma^{-1}$.\qed
	\end{lemma}
	In a similar manner we have the following.
	\begin{proposition}
		\label{prop01}
		Let $f=f_1f_2\cdots f_n\in \F_n, \phi(f)=\sigma$, and $i \in [n]$. 
		\begin{enumerate}
			\item $\Rmo(f)=\AX(\sigma)$.
			\item $f_i=\sigma(i)$ iff $i\in \AX(\sigma)$.
		\end{enumerate}
	\end{proposition}
	\begin{proof}
		\textbf{Point 1.} In the product $(n,f_n)\cdots(2,f_2)(1,f_1)$ the transpositions $(j,f_j)$ with $j<i$ fix the integer $i$; the transposition $(i,f_i)$ transforms $i$ into an integer $f_i\leq i$. If $i\notin\AX(\sigma)$, i.\ e., $\sigma(i)>i$, then the integer $f_i$ must appear in one of the terms $(j,f_j)$ for $j>i$. Therefore, there exists $j>i$ such that $f_i=f_j$, i.\ e., $i\notin\Rmo(f)$. Conversely, if $i\notin\Rmo(f)$, i.\ e., there exists $j>i$ such that $f_i=f_j$, then the transposition $(i,f_i)$ transforms $i$ into an integer $f_i$ and the transposition $(j,f_i)$ transforms $f_i$ into a $j$. Now $j$ can appear in the transpositions $(\ell,f_\ell)$ with $\ell>j$ if and only if $j=f_\ell\in\Image(f)$. In this case $j$ is transformed into a larger integer $\ell$, and so on. Thus, $\sigma(i)>i$ and $i\notin\AX(\sigma)$.
		
		\textbf{Point 2.} Suppose $\sigma$ has an anti-excedance at $i$. Then $\sigma(i)\leq i$. Observe that the procedure to compute $\phi^{-1}(\sigma)$ does not move $\sigma(i)$ during its first $n-i$ iterations. Indeed each of the steps for $k=n, n-1, \ldots, i+1$ swaps  the integers $k$ and $\sigma(k)$. None of the integers involved in the swaps can be  $\sigma(i)$, because on one hand $k$ only take values larger than $i$ and $\sigma(i) \leq i$ and on the other hand none of the $\sigma(k)$ can be equal to $\sigma(i)$ because $\sigma$ is a bijection. Thus $\sigma_i(i)=\sigma(i)$,  and hence $f_i=\sigma(i)$.
	\end{proof}
	\section{Transposition arrays  and cycles}
	\begin{definition} 
		Let $\sigma\in \mathfrak{S}_n$. For any integer $i\in [n]$ we define the \textit{nearest orbital minorant of $i$} (under $\sigma$) as the integer $j\leq i$ such that $j=\sigma^t(i)$ with $t\geq1$ chosen as small as possible. We denote this integer by $\nom_\sigma(i)$ or simply $\nom(i)$ when there is no ambiguity regarding the permutation $\sigma$.
	\end{definition}
	We now provide an alternative definition of the transposition array $\phi^{-1}(\sigma)$ that is related to the cycle structure of the permutation $\sigma$. It suffices to observe the effect of multiplying the initial  permutation by the transpositions $(i,\sigma(i))$.
	
	\begin{theorem}
		For every $\sigma=\sigma(1)\sigma(2)\cdots\sigma(n)\in \mathfrak{S}_n$ we have  $\phi^{-1}(\sigma)=f=f_1f_2\cdots f_n$, where
		$
		f_i=\nom(i), \forall i\in [n].$
		\label{prop03.13}
	\end{theorem}
	\begin{proof}
		We use descending recursion. In other terms, to prove that the claim is true for all $i \in [n]$, we show that:
		\begin{enumerate}
			\item it is true for $n$  and
			\item if it is true for all integers $j$ with $i < j \leq n$, then it is true for $i$.
		\end{enumerate}
		The claim is trivially true for $i=n$ since $f_n=\sigma(n)$ and $\sigma(n) \leq n$.
		Then let $i <n$. By Proposition \ref{prop01} (2) the claim is trivially true if $i$ is an anti-excedance (with $t=1$, where $t$ is as above).
		Then let us suppose that $\sigma$ has an excedance at $i$, say $\sigma(i)=j_1$ with $j_1>i$. Observe that during the iteration step (of the procedure to construct $f$ from $\sigma$) corresponding to $j_1$, the multiplication of the permutation by the  transposition $(j_1,\sigma_{j_1}(j_1))$ changes the image of $i$ from $j_1$ to $\sigma_{j_1}(j_1)= f_{j_1}$. But by recurrence hypotheses, $f_{j_1}=\sigma^{t_1}(j_1)$, where $t_1$ is the smallest integer such that $j_2:=\sigma^{t_1}(j_1)<j_1$.
		We can then define a sequence of integers $j_1 > j_2 > \cdots >j_p >  j_{p+1}$  such that :
		\begin{itemize}
			\item $ j_p >  i  \geq j_{p+1}$;
			\item $j_1  = \sigma(i)$  and $j_{p+1} = \sigma_{i}(i)= f_i$;
			\item $j_{k+1}=\sigma_{j_k}(j_k)=\sigma^{t_k}(j_k)$ for $1 \leq k \leq p$ and $t_k$ is the smallest positive integer such that $\sigma^{t_k}(j_k)< j_k$. (Thus $j_{k+1}=\sigma^{1+ t_1+ \cdots + t_k}(i)$.)
		\end{itemize}
		Therefore, if  $t= 1 + \sum_{k=1}^p t_k$ we have  $ f_i=j_{p+1} = \sigma^{t}(i)$, we only have to prove that $t$ is the smallest integer such that $\sigma^{t}(i) \leq i$.
		Suppose there exists an integer $s<t$ such that $\sigma^{s}(i) \leq i$. Hence there exists an integer $p_0 \leq p$ such that $$1 + \sum_{k=1}^{p_0-1}  t_k < s < 1 + \sum_{k=1}^{p_0} t_k.$$
		Let $r = s - \left(1 + \sum_{k=1}^{p_0 - 1}  t_k \right)$,
		then
		$$\sigma^{s}(i) =  \sigma^{\left(1 + \sum_{k=1}^{p_0 - 1}  t_k\right) + r }(i) = \sigma^r (\sigma^{1 + \sum_{k=1}^{p_0 -1}  t_k } (i))=\sigma^r(j_{p_0}),$$
		this integer is smaller than $i$, therefore it is smaller than $j_{{p_0}}$, but $r$ is smaller  than $t_{p_0}$ and this is a contradiction because, by recurrence hypothesis $t_{p_0}$ is the smallest integer such that $\sigma^{t_{p_0}}(j_{p_0}) \leq j_{p_0}$.
	\end{proof}
	
	When $\sigma$ is represented graphically as the reunion of cyclic graphs (each corresponding to a cycle), to find $\nom(i)$ it suffices to start from $i$ and follow the arcs of the cycle (connecting $i$ and $\sigma(i)$) until one meets an integer $j\leq i$. Obviously, when (and only when) $i$ is the minimum of its own cycle, one has $\nom(i)=i$.
	
	Theorem \ref{prop03.13} suggests also a new interpretation of the integer $f_i$ and a new algorithm to go from $f$ to $\sigma$ (as a product of disjoint cycles and not as a word).
	\begin{proposition}
		If $f\in \F_n$,  then $\sigma =\phi(f)$ can be constructed as follows.
		For $i=1,2,\ldots,n$:
		\begin{itemize}
			\item if $f_i = i$, then add a new singleton cycle: $(i)$
			\item if $f_i < i$, then insert $i$ before $f_i$ in its cycle.
		\end{itemize}
	\end{proposition}
	\begin{example}
		Take $f=1132532\in \F_7$. Then $\phi(f)=\sigma$ can be obtained as follows:
		\begin{center}   
			$    \begin{array}{l}
				(1)\\
				(2,1)\\
				(2,1)(3)\\
				(4,2,1)(3)\\
				(4,2,1)(3)(5)\\
				(4,2,1)(6,3)(5)\\
				(4,7,2,1)(6,3)(5)=\phi(f)=\sigma.
			\end{array}$
		\end{center} 
	\end{example}
	The above construction shows that the number of fixed points of $f$ and the number of cycles of $\sigma=\phi(f)$ are equal. 
	
	Theorem \ref{prop03.13} implies in particular that the integers $i$ and $f_i$ are always in the same cycle of the permutation for all $i, 1 \leq i \leq n$.
	\begin{corollary} Let $f=f_1f_2\cdots f_n\in\F_n$ and $\phi(f)=\sigma$.
		\begin{enumerate}
			\item $i\in \Fxdp(f)$ iff $i$ is a minimum of a cycle of $\sigma$.
			\item If $f_i=f_j$, then $i$ and $j$ are in the same cycle of $\sigma$.
		\end{enumerate}
	\end{corollary}	
	\section{Non-Decreasing transposition arrays} 
	We say that a subexcedant function $f=f_1f_2\cdots f_n$ is \emph{non-decreasing} if $1 = f_1\leq f_2\leq \cdots \leq f_n$. Let $F_n^{\scriptscriptstyle\nearrow}$ denote the set of all non-decreasing subexcedant functions over $[n]$. For $n\ge1$, we have \begin{eqnarray}\label{eqnCat}
		|F_n^{\scriptscriptstyle\nearrow}|=c_n=\frac{1}{n+1}\binom{2n}{n},
	\end{eqnarray}the Catalan number (the set $F_n^{\scriptscriptstyle\nearrow}$ and the set of \emph{N-E lattice paths} from the point $(0,0)$, the origin, to the point $(n,n)$ and never going above the diagonal are in bijection (see \cite[Lemma $1.32$]{Bo})). 
	
	We let $\mathfrak{S}_n^{\scriptscriptstyle\nearrow}:=\{\sigma\in\mathfrak{S}_n : \phi^{-1}(\sigma)\in\F_n^{\scriptscriptstyle\nearrow}\}$, i.\ e., the set of permutations whose transposition array is non-decreasing. For instance, $F_3^{\scriptscriptstyle\nearrow}=\{111, 112, 113, 122, 123\}$ and $\mathfrak{S}_3^{\scriptscriptstyle\nearrow}=\{231, 312, 213, 132, 123\}$.
	
	Here we present a characterization for the permutations having non-decreasing transposition arrays under the bijection $\phi^{-1}$ and describe a certain number of their properties.
	
	If $f=f_1f_2\cdots f_n$ is a non-decreasing subexcedant function, then we can encode $f$ by the vector $r=(r_1, r_2, \ldots, r_n)$, where $r_i$ is the number of occurrences of the letter $i$ in $f$. These vectors have the following characterization.
	\begin{eqnarray*} 
		0\leq r_i\leq n-i+1, \forall i\in[n], r_1\neq0~~\mbox{and}~\sum_{i=1}^nr_i=n.
		\label{eqn01}
	\end{eqnarray*}
	\begin{proposition}
		If $\sigma=\sigma(1)\cdots\sigma(n{-}1)\in\mathfrak{S}_{n-1}$ is associated to a non-decreasing transposition array, and $\sigma'=(n,j)\circ\sigma$, where $\sigma(n{-}1)\leq j\leq n$, then so is $\sigma'$.
	\end{proposition}
	\begin{proof}
		If $f=f_1f_2\cdots f_{n-1}$ is the transposition array of $\sigma$, then the transposition array of $\sigma'$ is $f'=f_1f_2\cdots f_{n-1}j$. Since $f_{n-1}=\sigma(n{-}1)\le j$, $f'$ is non-decreasing.
	\end{proof}
	\begin{proposition}
		\label{CharNonDecPerm}
		Let $\sigma\in\mathfrak{S}_n$ and $f$ be its transposition array. Then, $f$ is non-decreasing if and only if
		\begin{enumerate}
			\item the subword of anti-excedance letters of $\sigma$ is increasing. In other terms, for any two anti-excedances $i$ and $j$ of $\sigma$ with $i<j$, one has  $\sigma(i)<\sigma(j)$, and
			\item  elements of $[n]$ having the same nearest orbital minorant  under $\sigma$ form an integer interval, i.\ e., if $j\in\Image(f)$, then $f^{-1}(j)=[a,b]=\{a,a{+}1,\ldots, b\}$ for some integers $a\le b$.
		\end{enumerate}  
	\end{proposition}
	
	\begin{proof}
		Let $r=(r_1, r_2, \ldots, r_n)$ be the vector code of $f$. By \cite[Proposition 3.5]{Ma-Ra}, the elements of $\Image(f)$ are the  anti-excedances of $\sigma$, and these are the integers $k$ such that $r_k \neq 0$. Furthermore, by Proposition \ref{prop01}, only the rightmost occurrences of each element of $\Image(f)$ in $f$ correspond to the (positions of the) anti-excedances. If $f$ is non-decreasing, then for each $k$ such that $r_k \neq 0$, that position is exactly $\sum_{\ell=1}^k r_k$. Further, observe that the nearest orbital minorants of the elements of $[n]$ under $\sigma$ are elements of $\Image(f)$. If $i<k<j$, where $i, j, k\in[n]$ with $\nom(i)=\nom(j)$, then $\nom(i)= \nom(k)$ since $f$ is non-decreasing.
		
		Conversely, it is obvious that the subexcedant function $f=\phi^{-1}(\sigma)$ with $\sigma$ having the given conditions is non-decreasing.
	\end{proof}	
	\begin{proposition}
		If $\sigma\in\mathfrak{S}_n^{\scriptscriptstyle\nearrow}$, then the anti-excedances of $\sigma$ and their images entirely characterize $\sigma$.
		\label{charNonDec}
	\end{proposition}
	
	\begin{proof}
		Let $\{i_1 < i_2 < \cdots < i_k{=}n\}$ be the anti-excedances of $\sigma$ and $\{j_1{=}1 < j_2 < \cdots < j_k\}$ their respective images under $\sigma$ (with $i_t \geq j_t, \forall t\in[k]$).	By Proposition  \ref{prop01}, $i_t$ is the rightmost occurrence such that $f_{i_t}=j_t$ and the number of positions having $j_t$ as a value is $i_t {-} i_{t-1}$ (if one poses $i_0=0$).  Therefore, $f$ is entirely determined and consequently $\sigma$ is too. \end{proof}
	\begin{example}
		Let $\sigma \in \mathfrak{S}_9^{\scriptscriptstyle\nearrow}$  with $\AX(\sigma)=\{3, 4, 6, 9\}$ and $\AXL(\sigma)=\{1, 2, 3, 5\}$, respectively. If $f= \phi^{-1}(\sigma)$, then we have $f=11\mathbf{12}3\mathbf{3}55\mathbf{5}\in \F_9^{\scriptscriptstyle\nearrow}$, and thus $\sigma= 471263895$.
	\end{example}	
	Let $T_n$ be the set of tuples of integer pairs  $\langle (i_1,j_1), (i_2,j_2), \ldots, (i_k,j_k)\rangle$ satisfying the following:
	
	\begin{align*}\label{Catalantuples}
		(i)~~& 1 \leq k \leq n,\notag\\
		(ii)~~&1 \leq i_1 < i_2 < \cdots < i_k=n,\notag\\
		(iii)~~& 1 = j_1 < j_2 < \cdots < j_k, \text{ and } \\
		(iv) ~~&j_s\leq i_{s-1}+1 \text{ for }2\leq s\leq k.\notag
	\end{align*}
	We can deduce the following corollary as a new proof of a result already known as an incarnation of Catalan numbers (see \cite{Hu,St1}).
	\begin{corollary}
		For $n\ge0$, we have $|T_n|=c_n$. \end{corollary}
	\begin{proof}
		Two such sequences $\langle i_1, i_2, \ldots i_k\rangle$  and $\langle j_1, j_2, \ldots j_k\rangle$ constitute respectively the set of anti-excedances and their images of a permutation associated  with a non-decreasing transposition array, namely $ 1^{i_1} j_2^{i_2-i_1} \cdots  j_k^{n-i_{k-1}}$, and by (\ref{eqnCat}) we have the result.
	\end{proof}
	\section{Avoidance of length 3 patterns in  $\mathfrak{S}_n^{\scriptscriptstyle\nearrow}$}
	In this section we study length $3$ pattern avoidance in the set $\mathfrak{S}_n^{\scriptscriptstyle\nearrow}$. We denote the set of permutations in $\mathfrak{S}_n^{\scriptscriptstyle\nearrow}$  avoiding a pattern $\pi$ by $\mathfrak{S}_n^{\scriptscriptstyle\nearrow}(\pi)$. 
	\subsection*{The case $123$.} The first few terms of the sequence of the numbers of $123$-avoiding permutations over $[n], n\geq1$, corresponding to non-decreasing transposition arrays is $1, 2, 4, 4, 3, 0, 0, 0, \ldots$.
	
	\begin{proposition}
		We have $|\mathfrak{S}_n^{\scriptscriptstyle\nearrow}(123)|=0$ for $n\geq6$.
	\end{proposition}
	\begin{proof}
		If $\sigma\in\mathfrak{S}_n^{\scriptscriptstyle\nearrow}$ such that $\ax(\sigma)\geq3$, then $\sigma$ must contain the pattern $123$ since the subword of anti-excedance letters of $\sigma$ is increasing. Let $\ax(\sigma)\leq 2$. If $\ax(\sigma){=}1$, then $\sigma=234\cdots n1$ and it contains the pattern $123$ for all $n\ge4$. Now suppose that $\ax(\sigma){=} 2$ and let $\AX(\sigma)=\{i, n\}$ such that $i<n$. Then $\sigma(i){=}1<\sigma(n)$. Let $\sigma(n)=j\ge2$. Then, the transposition array of $\sigma$ has the form $f=11\cdots1jj\cdots j\in\F_n^{\scriptscriptstyle\nearrow}$. If $j>i$, then $j=i{+}1$ and thus $\sigma=\phi(f)=23\cdots i1(i{+}2)(i{+}3)\cdots n(i{+}1)$. So removing $1$ and $i{+}1$ leaves an increasing sequence of length $n{-}2$. Similarly, if $j\le i$, then $\sigma=23\cdots(j{-}1)(i{+1})(j{+}1)\cdots i1(i{+}2)(i{+}3)\cdots nj$ and removing $i{+}1, 1$ and $j$ leaves an increasing sequence of length $n{-}3$, and the pattern $123$ occurs in $\sigma$ for $n\ge6$. 
	\end{proof}
	\subsection*{The case $132$.} In a subexcedant function $f$, a plateau $i$ of $f$ is said to be of height $h$ if $f_i=f_{i+1}=h$. Further, $f$ has a \emph{stair} at $i\in[n{-}1]$ if $f_{i+1}=f_i+1\ge3$. A \emph{staircase} in $f$ is the letters of a maximal interval of integers in which all the integers except the last are stairs.
	
	Let $V_n:=\{f\in F_n^{\scriptscriptstyle\nearrow} : f \text{ may have plateaus only of height }1\}$, i.\ e., $f \in V_n$ iff $f_1=f_2=\cdots=f_k=1<f_{k+1}<f_{k+2}<\cdots<f_n$, where $1\leq k\leq n$. For instance, $V_3= \{111, 112, 113, 123\}$.
	\begin{proposition}
		If $\sigma\in\mathfrak{S}_n^{\scriptscriptstyle\nearrow}(132)$, then $f=\phi^{-1}(\sigma)\in V_n$.
	\end{proposition}
	\begin{proof}
		Let $\sigma\in\mathfrak{S}_n^{\scriptscriptstyle\nearrow}(132)$ and suppose that $f = \phi^{-1}(\sigma) \notin V_n$, i.\ e., $f$ has plateau(s) of height $h$ with  $h \geq 2$. Let $k$ be the rightmost position such that $f_k=1$ and $p>k+1$ be the rightmost position such that $f_p=h$. By proposition \ref{prop01}, $\sigma$ has an anti-excedance in $k$ with image value $1$ and an anti-excedance in $p$ with image value $h\leq p$. In the position $p-1$, we have $f_{p-1}=h$ and an excedance for $\sigma$, i.\ e., a letter  larger than $p-1$, therefore $\sigma$ would have a pattern $132$ in the positions $k, p{-}1$ and $p$, a contradiction. 
	\end{proof}
	\begin{remark}
		Let $f \in V_n$, where $\plat(f)=k$ and $\sigma=\phi(f)$. Then
		\begin{enumerate}
			\item for all $i$ with $k{+}1\leq i\leq n$, we have $\sigma(i)=f_i$, because all these are anti-excedance positions. Thus, $\AXL(\sigma)$ is an increasing sequence of integers, and $\Exc(\sigma) = [k]$;
			\item if there exists $i>1$ such that $f_i=i$, then $f_j=j, \forall j>i$, in particular $f_n=n$.
		\end{enumerate}
	\end{remark} 
	\begin{definition}\label{defdecompose}
		Let $f\in V_n$, where $\plat(f)=k$. We decompose $f=1| H_0| L_1| \cdots| L_p$, where $H_0=1^{k_0}$, $k_0=k$, and the $L_i$'s are staircases. If $H_i$ is the last block of $L_{j-1}$, then we cut $L_j$ into blocks of size $k_i$ (with a possible remainder of size $(\ell_j\mod k_i)$ and let $k_{i+1}=((\ell_j-1)\mod k_i)+1$, where $\ell_j=|L_j|$. Then, let the resulting decomposition be $f=1| H_0| H_1|H_2|\cdots|H_s$ with $H_i$ having size $k_i, 0\le i\le s$. 
	\end{definition}
	For instance, take $f=1~1~1~1~2~3~4~5~6~7~8~9~11~12$. We have $L_0=2~3~4~5~6~7~8~9, L_2=11~12, k=k_0=3, \ell_1=8$, and $\ell_2=2$. First, we cut $L_1$ into two blocks of size $k_0=3$ and the block of size the remainder $(5\mod3)=2$. So $k_1=k_2=3, k_3=2$. Next, we cut $L_2$ into block of size $k_3=2$ since the remainder $(2\mod2)=0$. So $k_4=2$. Thus, $f$ becomes $f=1| H_0|H_1|H_2|H_3|H_4=1|1~1~1|2~3~4|5~6~7|8~9|11~12$ (see the following figure). 
	\begin{figure}[ht!]
		\centering
		\includegraphics[width=4.1cm]{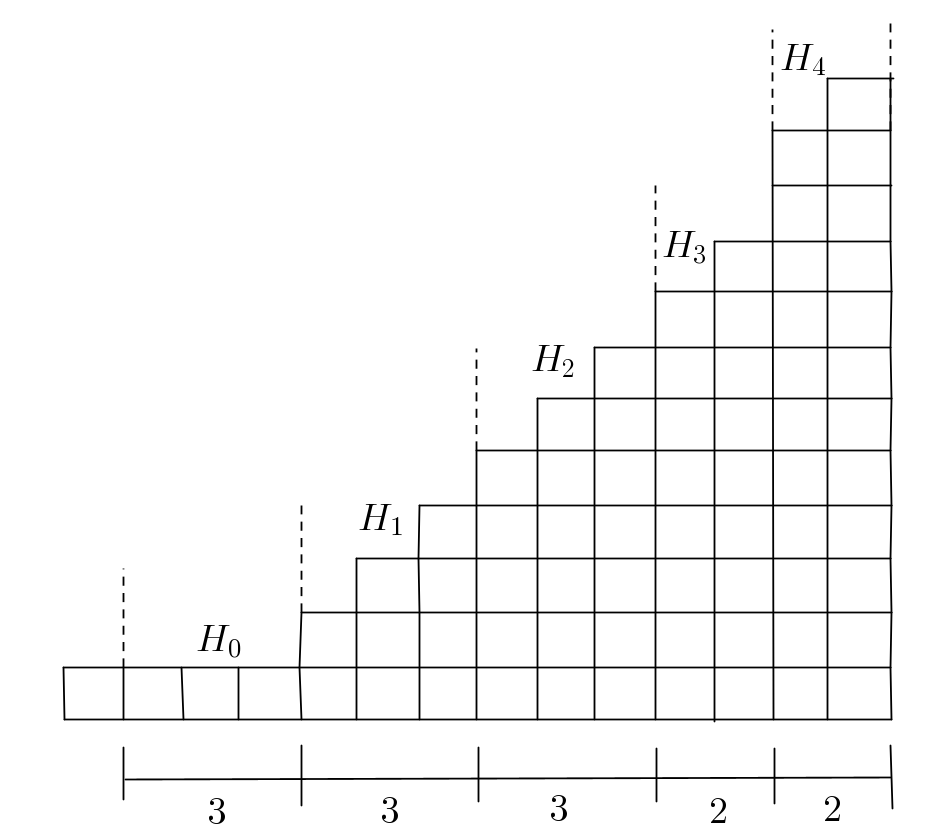}
	\end{figure}
	
	Let $\mathfrak{D}_n:=\{f\in V_n :\sigma=\phi(f) \in\mathfrak{S}_n^{\scriptscriptstyle\nearrow}(132), \Fxdp(f)=\{ 1\}\}$, and $\mathfrak{D}_{n,k}:=\{f\in\mathfrak{D}_n: \plat(f)=k\}$. Further, let $d_n:=|\mathfrak{D}_n|$ and $d_{n,k}:=|\mathfrak{D}_{n,k}|$. The following proposition characterizes the elements of $\mathfrak{D}_n$.
	\begin{proposition}
		\label{propChar132sef}
		Let $f\in V_n$ be given in its decomposed form $f=1| H_0| H_1|H_2|\cdots|H_s$ as in Definition \ref{defdecompose}.
		Then, $f\in\mathfrak{D}_n$ if and only if
		\begin{eqnarray}\label{eqnsef132avoid}
			f_{k_0+k_1+\cdots+k_i+2}= k_0+k_1+\cdots+k_{i-1}+2,~ 1\leq i\leq s,\end{eqnarray}
		i.\ e.,  each block $H_i$ with $i>0$ begins with the integer equal to the initial position of the block $H_{i-1}$.
	\end{proposition}
	\begin{proof}
		Suppose $f\in\mathfrak{D}_n$, where $\plat(f)=k_0$.
		We shall use induction to prove that $f$ satisfy (\ref{eqnsef132avoid}) and that 
		\begin{eqnarray}\label{eqnsigma}
			\sigma_{1+k_0+\cdots+k_i}(j)=1+k_0+\cdots+k_{i-1}+j, \end{eqnarray}
		where $j=1, 2, \ldots, k_i$, and ($\sigma_a\in\mathfrak{S}_a, 1\le a\le n$, as given in Section \ref{insertionmthd}).
		
		Observe that $\sigma_{1+k_0}=\phi(1H_0)=23\cdots(k_0{+}1)1$, and thus (\ref{eqnsigma}) holds for $i=0$. We show that $f_{k_0+2}=2$. If $f_{k_0+2}>2$, then the integer $\sigma_{1+k_0}(1){=}2$ remains unmoved in the successive steps of the construction of $\sigma$ because $2\notin\Image(f)$. So $\sigma(1)=2$. Since $f_n<n$ and $\Exc(\sigma)=[k_0]$, we have some $j\in[2,k_0]$ such that $\sigma(j)=n$. This implies that $\sigma$ contains the pattern $132$ in the positions $1, j, k_0+2$ and this is a contradiction. Thus, instead, $f_{k_0+2}=2$ and (\ref{eqnsef132avoid}) holds for $i=0$.
		
		Now assume that (\ref{eqnsef132avoid}) and (\ref{eqnsigma}) hold for $i=q{-}1$. We prove it for $i=q$. If we let $t:=k_0+k_1+\cdots+k_q+2$ (the inial position of the block $H_{q+1}$), then $k_0+k_1+\cdots+k_{q-1}+2=t{-}k_q$.	By (\ref{eqnsigma}) (for $i=q{-}1$) we have the largest $k_{q-1}$ integers already inserted into the first $k_{q-1}$ positions of $\sigma_{t-k_q-1}$ in increasing order. Since $H_q=[t{-}k_q{-}k_{q-1},t{-}k_{q-1}{-}1]$ (by (\ref{eqnsef132avoid})) and $k_{q-1}\geq k_q$, the next $k_q$ steps insert the largest $k_q$ integers in $\sigma_{t-1}$ in increasing order. Therefore, (\ref{eqnsigma}) holds. 
		
		Now we prove (\ref{eqnsef132avoid}), i.\ e., $f_t=t{-}k_q$. Since $f\in V_n$, we have $f_t\ge f_{t-1}{+}1$. If $f_t=f_{t-1}{+}1$, then $t{-}1$ is a stair, whence $k_{q-1}=k_q$, and hence $f_t=f_{t-k_q}+k_q=(t-k_q-k_{q-1})+k_{q-1}=t{-}k_q$ indeed. Thus, instead assume that $f_t> f_{t-1}{+}1$. By (\ref{eqnsigma}) we have $\sigma_{t-k_q-1}(j)=t{-}k_q{-}k_{q-1}{-}1{+}j, 1\le j\le k_{q-1}$. If $f_t<t{-}k_q$, then there is some $j\in[k_q{+}2,k_{q-1}]$ such that $\sigma_{t-1}(j)=\sigma_{t-k_q-1}(j)=f_t$. Since $\sigma_t$ is obtained from $\sigma_{t-1}$ by applying the transposition $(t, f_t)$, we have $\sigma_t(j)=t$ and $\sigma_t(t)=f_t$. Thus, $\sigma_t(k_q{+}1)=\sigma_{t-k_q-1}(k_q+1)=t{-}k_{q-1}=f_{t-1}{+}1<f_t=\sigma_t(t)<\sigma_t(j)$ and $k_q{+}1<j<t$. Since in the successive steps $\sigma_t(k_q{+}1)$ and $\sigma_t(t)$ remain unmoved and $\sigma_t(j)$ can only be replaced by a larger integer, we have the pattern $132$ in $\sigma$ and this is a contradiction.
		If $f_t>t{-}k_q$, then $f_t$ is one of the largest $k_q{-}1$ integers already inserted in positions $[2,k_q]$ and the positions $1, \sigma_{t-1}^{-1}(f_t)$ and $t$ would create the pattern $132$. Therefore, $f_t=t{-}k_q$.
		
		Conversely, suppose that $f$ satisfies the given condition. We can decompose the block of excedance letters of $\sigma$ (which are the letters in the first $k_0$ positions) into blocks $P_1|P_2|\cdots|P_u$ made of increasing numbers and such that for every $a\in P_j, b\in P_{j+1}$ we have $a>b$, where $j=1, 2, \ldots, u-1$. Indeed for all $r$ with $1\leq  r \leq s$, it can be verified that if $\tau= \phi(1| H_0| H_1|H_2|\cdots|H_r)$, then $\tau$ is obtained from $\tau'= \phi(1| H_0| H_1|H_2|\cdots|H_{r-1})$ by replacing the integers $\tau'(1), \ldots, \tau'(k_r)$ with the ${k_r}$ largest integers in increasing order and appending $\tau'(1), \ldots, \tau'(k_r)$ at the end. Thus, if $\sigma$ had an occurrence $132$, then the $3$ and the $2$ could not both be among the anti-excedance letters because they form an increasing sequence (not necessarily an interval). Hence at least the $1$ and the $3$ must be in the excedance block, and thus in the same block $P_j$ for some $j$ because the $P_j$'s are intervals arranged in decreasing order of their first elements. So there can be no integer playing the role of $2$ outside $P_j$. Thus there is no such occurrence.
	\end{proof}
	\begin{theorem}
		\label{thm132avoid}
		The number $a_n:=|\mathfrak{S}_n^{\scriptscriptstyle\nearrow}(132)|$ satisfies the recurrence relation:
		\begin{eqnarray} \label{ReccurrencePartition}a_n=a_{n-1}+p_{n-1}, n\geq1,~~ a_0=0, p_0=1,\end{eqnarray}
		where $p_n$ is the number of integer partitions of $n$.
	\end{theorem}
	\begin{proof}
		Note that $a_n$ is the number of non-decreasing transposition arrays $f$ over $[n]$ such that $\sigma=\phi(f)\in\mathfrak{S}_n^{\scriptscriptstyle\nearrow}(132)$.	
		The term $a_{n-1}$  in the right-hand side of the recurrence relation can be explained easily. Let $f=f_1f_2\cdots f_{n-1}$ such that $\sigma=\phi(f)\in\mathfrak{S}_{n-1}^{\scriptscriptstyle\nearrow}(132)$. Then the function $f'$ obtained from $f$ by appending $n$ at its end corresponds to a permutation $\sigma'=\phi(f')\in\mathfrak{S}_n^{\scriptscriptstyle\nearrow}(132)$. This term counts the non-decreasing subexcedant functions with fixed points and such that $\phi(f')\in\mathfrak{S}_n^{\scriptscriptstyle\nearrow}(132)$  and it contributes $a_{n-1}$ to $a_n$. 
		For the  term $p_{n-1}$, we show in the following lemmas and propositions that the number of functions $f$ with no fixed points greater than $1$, with $\sigma=\phi(f)\in\mathfrak{S}_n^{\scriptscriptstyle\nearrow}(132)$, equals the number $p_{n-1}$ of all integer partitions of  $n-1$ and more precisely that the number $d_{n,k}$ equals the number of integer partitions over $[n-1]$ with the largest part equal to $k$ (see \cite{St1}). Thus, this contributes $p_{n-1}$ to the number $a_n$.
	\end{proof}
	\begin{remark}
		The sequence of the numbers $a_n$ is the same as the sequence given in OEIS number \underline{A000070} but shifted by one step. In particular, it has
		has an ordinary generating function $$\prod_{j\ge1}\frac{x}{(1-x^j)^{1+\delta_{1,j}}}.$$
	\end{remark}
	\begin{lemma} 
		\label{lemmaPermConst1}
		Let $f\in\mathfrak{D}_{n,k}$ and $f'$ be obtained from $f$ by increasing each integer greater than or equal to $k+2$ by $1$ and inserting a $1$ at the beginning, then  $f'\in\mathfrak{D}_{n+1,k+1}$. 
	\end{lemma}
	\begin{proof}
		If $f$ is decomposed as in Definition \ref{defdecompose} we have $f_{k_0+k_1+2}=k_0+2$. If $f'=f_1'f_2'\cdots f_{n+1}'$, then 
		$$f_j'=\begin{cases}
			1, &\mbox{ if } 1\leq j\leq k_0+2;\\
			f_{j-1}, &\mbox{ if } k_0+3\leq j\leq k_0+k_1+2;\\
			f_{j-1}+1, &\mbox{ if } k_0+k_1+3\leq j\leq n+1.
		\end{cases}$$
		Let us show that $f'\in\mathfrak{D}_{n+1,k+1}$, i.\ e., $f'$ satisfies the condition in Proposition \ref{propChar132sef}. The number of blocks of $f$ and $f'$ are the same but the size of $H_0'$ of $f'$ is increased by $1$, where $f'=1| H_0'| H_1'| H_2'|\cdots| H_s'$. So $f_{k_0'+2}'=f_{k_0+3}'=f_{k_0+2}=2$, and $f'_{k_0'+k_1'+\cdots+ k_i'+2}=f_{k_0+1+k_1+\cdots +k_i+2}+1=k_0+k_1+\cdots+ k_{i-1}+3=k_0'+k_1'+\cdots+k_i'+2, i\geq 1$. Thus, $f'\in\mathfrak{D}_{n+1,k+1}$.
	\end{proof}
	\begin{lemma} 
		\label{lemmaPermConst2}
		Let $f\in\mathfrak{D}_{n,k}$ and $f'$ be obtained from $f$ by increasing each integer greater than $1$ by $k$ and then inserting the subword $23\cdots k+1$ after the rightmost occurrence of $1$. Then  $f'\in\mathfrak{D}_{n+k,k}$. 
	\end{lemma}
	\begin{proof}  
		Since $((\ell_1-1)\mod k)+1=((\ell_1+k-1)\mod k)+1$ the insertion of the subword $23\cdots k+1$ in the operation creates a new block of size $k$ in $f'$ after $H_0$. That is, $f'=1| H_0| H_1'| H_2'|\cdots| H_{s+1}'$, where $H_1'=2~3~4~\cdots~k+1$ and $a+k\in H_i'$ if $a\in H_{i-1}, i\geq2$. So, $f_{k_0+2}'=2$ and \begin{align*}
			f'_{k_0'+k_1'+\cdots +k_i'+2}=&f_{k_0+k_0+k_1+\cdots+ k_{i-1}+2}'\\
			=&f_{k_0+k_0+k_1+\cdots+ k_{i-1}+2-k_0}+k_0\\
			=&(k_0+k_1+\cdots+ k_{i-2}+2)+k_0\\
			=&k_0'+k_1'+\cdots +k_{i-1}'+2, ~~~i\geq 1.\end{align*} Thus, $f'\in\mathfrak{D}_{n+k,k}$. 
	\end{proof}
	From Lemma \ref{lemmaPermConst1} and Lemma \ref{lemmaPermConst2} we deduce that each function $f$ in $\mathfrak{D}_{n,k}$, can uniquely be obtained either from a function in $\mathfrak{D}_{n-1,k-1}$  by the construction in Lemma \ref{lemmaPermConst1} or  from a function in $\mathfrak{D}_{n-k,k}$  by the construction in Lemma \ref{lemmaPermConst2}. Furthermore, the two sets obtained this way are disjoint because those obtained by the construction in Lemma \ref{lemmaPermConst1} do not have $k+2$ in their image, while those obtained by the construction in Lemma \ref{lemmaPermConst2} do. As a result, we have the following proposition.
	\begin{proposition}
		\label{proprecsefplat}
		The number $d_{n,k}$ satisfies the recurrence relation:
		\begin{eqnarray}
			d_{n,k}=d_{n-1,k-1}+d_{n-k,k},
		\end{eqnarray}
		where $d_{0,0}=1$ and $d_{n,k}=0$ if $n\le0$ or $k\le0$ and $n$ and $k$ are not both zero.\qed
	\end{proposition}
	This recurrence is the same as the one satisfied by the number of integer partitions of $n$ whose largest part is $k$ (see \cite[Section $1.7$]{St1}), we then deduce the following.
	\begin{corollary}
		For $n\ge1$, we have $d_{n}=p_{n-1}$ .\qed
	\end{corollary}
	\begin{remark}
		The number $d_{n,k}$ equals the number of permutations $\sigma\in\mathfrak{S}_n^{\scriptscriptstyle\nearrow}(132)$, where $\sigma(n)<n$, having $k$ excedances.
	\end{remark}
	It is possible to define a direct bijection between the set $\mathfrak{D}_n$ and the set   $P_{n-1}$ of all integer partitions of $n{-}1$. 
	Let $f\in\mathfrak{D}_n$ and $\rho: \mathfrak{D}_n\mapsto P_{n-1}$ be the map which associates $f$ with an integer partition $\rho(f)=\lambda\vdash n{-}1$ defined as follows.
	Decompose $f$ as indicated in Definition \ref{defdecompose}: $f=1| H_0| H_1| \cdots| H_s$, and call $k_i$ the size of $H_i$, then we set $\rho(f)=\lambda=k_0k_1\cdots k_s$, this is clearly a partition of $n{-}1$ with the largest part equal to $k_0$. 
	\begin{proposition}
		The map $\rho: \mathfrak{D}_n\mapsto P_{n-1}$ is a bijection.
	\end{proposition}
	\begin{proof} 
		We will show that $\rho$ can be inversed. Let $\lambda=\lambda_1\lambda_2\cdots\lambda_\ell$ be a partition of $n{-}1$. Let us define $\rho^{-1}(\lambda)=f$ over $[n]$ as follows. First start with $k=\lambda_1$ plateaus of height $1$, i.\ e., set $f_2=\cdots=f_{k+1}=1$. Then for $i\geq2$ insert $\lambda_i$ consecutive integers starting from the integer $2 + \sum_{j=1}^{i-2}\lambda_j $ (where, as usual $\sum_{j=1}^{0}\lambda_j=0$). 
	\end{proof}
	Let $\lambda=55533221$. We have $\lambda\vdash 26$. Starting from $k=\lambda_1=5$ we obtain
	$$f=1|1~1~1~1~1|2~3~4~5~6|7~8~9~10~11|12~13~14|17~18~19|20~21|23~24|25\in\mathfrak{D}_{27}.$$

	\subsection*{The case $213$.}
	We will prove that $|\mathfrak{S}_n^{\scriptscriptstyle\nearrow}(213)|=|\mathfrak{S}_n^{\scriptscriptstyle\nearrow}(132)|, n\geq 1$. 
	In order to do so we introduce an involution in the set $F_n^{\scriptscriptstyle\nearrow}$, called ``$\flip$''. This involution  associates each $f=f_1f_2\cdots f_n\in F_n^{\scriptscriptstyle\nearrow}$ with the function $\flip(f)=f_1'f_2'\cdots f_n'$ given as follows: $f_i'=n{+}1{-}\sum_{k=1}^{n-i+1} r_k$, where $(r_1, r_2, \ldots, r_n)$ is the $r$-vector of $f$. For instance, let $f=1~1~1~1~1~1~2~3~4~6~6~9~10~11~13$. The $r$-vector of $f$ is $(6, 1, 1, 1, 0, 2, 0, 0, 1, 1, 1, 0, 1, 0, 0)$, and hence $\flip(f)=1~1~1~2~2~3~4~5~5~5~7~7~8~9~10$. 
	\begin{figure}[ht!]
		\centering
		\includegraphics[width=9.5cm]{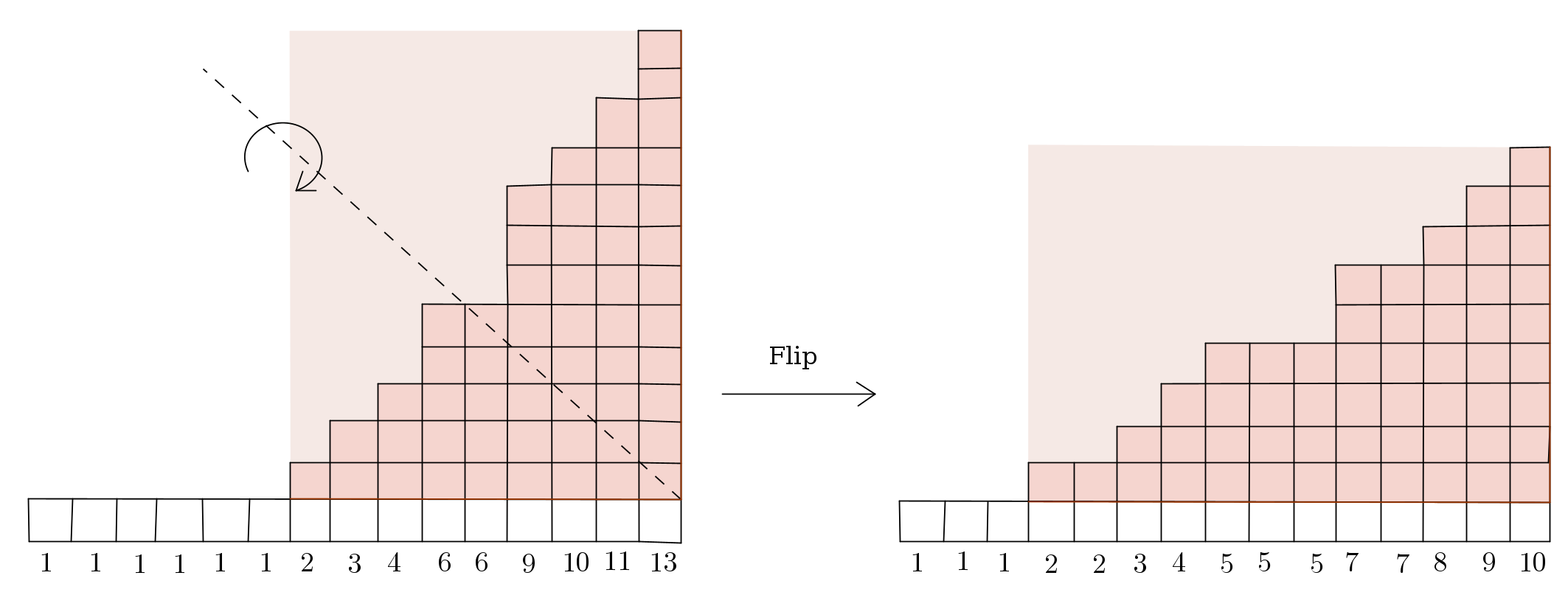}
	\end{figure}
	
	Although non-decreasing subexcedant functions can be seen as partitions of the integer $\sum_{i=1}^{n} f_i$ and although this operation is similar to the conjugation  in the set of integer partitions, the two operations are  different.
	\begin{remark}\label{remaplat}
		If $f\in \F_n^{\scriptscriptstyle\nearrow}$, then $\flip(f)$ starts with $n-f_n$ plateaus of height $1$.
	\end{remark}	
	\begin{lemma}\label{AddjandFlip}
		If $f'=f_1f_2\cdots f_{n-1}\in \F_{n-1}^{\scriptscriptstyle\nearrow}$ and $f=f_1f_2\cdots f_{n-1}j$, where $j\geq f_{n-1}$, then $\flip(f)$ is obtained from $\flip(f')$ by adding $1$ to the integers in the $j{-}1$ largest positions and inserting a $1$ at the beginning.
	\end{lemma}
	\begin{proof}
		Suppose that $f'$ has $r$-vector $(r_1, r_2, \ldots, r_{n-1})$ and $\flip(f')=g'=g_1'\cdots g_{n-1}'$. Then the $r$-vector of $f=f'\cdot j$ is $(r_1^*, \ldots, r_n^*)$, where $r_i^*=r_i{+}\delta_{i,j}, r_n^*=\delta_{n,j}$. So, we have $\flip(f)=g_1g_2\cdots g_n$, where $g_i=n+1-\sum_{k=1}^{n-i+1}r_k^*$. It is easy to see that $g_1=1, g_i=g_{i-1}'$ for $2\le i< n{-}j{+}2$, and $g_i=g_{i-1}'+1$ for $i\ge n{-}j{+}2$, i.\ e., $\flip(f)$ is exactly the function obtained from $\flip(f')$ by adding $1$ to the integers in the $j{-}1$ largest positions and inserting a $1$ at the beginning.
	\end{proof}
	\begin{theorem}
		If $f\in \F_n^{\scriptscriptstyle\nearrow}$ is the transposition array of the permutation $\sigma$, then $\flip(f)$ is the transposition array of the permutation $ \tau =((\sigma^{-1})^r)^c$, where $r$ and $c$ denote the operations of reverse and complement of a permutation respectively. 
	\end{theorem}
	\begin{proof}
		For any permutation $\pi \in \mathfrak{S}_n$ one has  $(\pi^r)^c=(\pi^c)^r= \psi_{n}\circ\pi\circ\psi_{n}$, where $\psi_{n}$ is the permutation $\left(\begin{array}{cccc}1&2&\cdots&n\\
			n&n-1&\cdots& 1\end{array}\right)$, therefore the equality to prove becomes: $\tau= \psi_{n}\circ \sigma^{-1}\circ\psi_{n}$.
		We prove the result by induction on $n$. The result is trivially true for $n=1$. Suppose $n > 1$ and that the result is true for $n{-}1$. Let $f=f' \cdot j$, where $f'$ is the prefix of length $n{-}1$ of $f$ and $j=f_n\geq f_{n-1}= f_{n-1}'$. Let $\sigma'$ be the permutation whose transposition array is $f'$ and let $\tau'$ be the permutation whose transposition array is $\flip(f')$. By induction hypotheses $\tau'= \psi_{n-1}\circ (\sigma')^{-1}\circ\psi_{n-1}$  or equivalently (since $\psi_{n-1}$ is an involution) $\psi_{n-1} \circ\tau' \circ\psi_{n-1} = (\sigma')^{-1}$. On the other hand, by the insertion method, $\sigma= (n, j)\circ\sigma'$ and hence $\sigma^{-1}= (\sigma')^{-1}\circ(n, j)$, therefore the relation to be proven becomes:
		\begin{eqnarray} 
			\tau &= & \psi_{n}\circ \psi_{n-1}\circ  \tau'\circ \psi_{n-1}\circ(n, j)\circ\psi_{n}.
			\label{flip}
		\end{eqnarray}  
		By the definition of the bijection $\phi$ we have
		$\tau'= \prod_{i=1}^{n-1} (n{-}i, f'_{n-i})$, then by Lemma \ref{AddjandFlip}, $\tau=\prod_{i=1}^{j-1} (n{-}i{+}1, f'_{n-i}+1)\circ \prod_{i=j}^{n-1} (n{-}i{+}1, f'_{n-i})$. By plugging these values in Equation \ref{flip} we obtain
		$$\prod_{i=1}^{j-1} (n{-}i{+}1, f'_{n-i}+1)\circ \prod_{i=j}^{n-1} (n{-}i{+}1, f'_{n-i})=\psi_{n}\circ \psi_{n-1}\circ  \left(\prod_{i=1}^{n-1} (n{-}i, f'_{n-i})\right)\circ \psi_{n-1}\circ(n, j)\circ\psi_{n}.$$
		In the symmetric group, if a permutation $\alpha$ is written as product of cycles and $\beta$ is another permutation, the conjugate $\beta \alpha \beta^{-1}$ can be computed by replacing every integer $i$ with $\beta(i)$ in the product of cycles giving $\alpha$.
		\begin{eqnarray*}
			\prod_{i=1}^{j-1} (n{-}i{+}1, f'_{n-i}+1)\circ \prod_{i=j}^{n-1} (n{-}i{+}1, f'_{n-i}) & = & \psi_{n}\circ \left(\prod_{i=1}^{n-1} (i, n{-}f'_{n-i})\right)\circ (n, j)\circ\psi_{n}\\
			& =& \prod_{i=1}^{n-1} (n{-}i{+}1, f'_{n-i}{+}1)\circ(1, n{-}j{+})\\
		\end{eqnarray*}
		
		by simplifying we get
		\begin{eqnarray*}
			\prod_{i=j}^{n-1} (n{-}i{+}1, f'_{n-i}) & = & \prod_{i=j}^{n-1} (n{-}i{+}1, f'_{n-i}{+}1)\circ(1, n{-}j{+}1)\\
		\end{eqnarray*}
		because of Remark \ref{remaplat}, all the $f'_i$ for $i=1, \ldots, n{-}j{+}1$ are equal to $1$.
		$$(n{-}j{+}1, 1)\cdots (3, 1)(2, 1)  = (n{-}j{+}1, 2) \cdots (3, 2)(2, 2)(1, n{-}j{+}1).$$
		It is straightforward to check that these two permutations are equal (and equal to the cycle $(1, 2,  \ldots, n{-}j{+}1)$).
	\end{proof}
	\begin{corollary}
		For $n\geq 1$, we have $|\mathfrak{S}_n^{\scriptscriptstyle\nearrow}(213)|=|\mathfrak{S}_n^{\scriptscriptstyle\nearrow}(132)|$.
	\end{corollary}
	\begin{proof}
		A permutation $\sigma$ contains the pattern $132$ if and only if $((\sigma^{-1})^r)^c$ contains the pattern $(((132)^{-1})^r)^c=213$, therefore $(\phi \circ \flip \circ ~\phi^{-1})$ is a bijection between $\mathfrak{S}_n^{\scriptscriptstyle\nearrow}(213)$ and $\mathfrak{S}_n^{\scriptscriptstyle\nearrow}(132)$.
	\end{proof}
	\subsection*{The case $231$.} 
	We do not have an explicit expression for the number of permutations of this class, but we can provide a lower bound showing that it grows exponentially with $n$. The first few terms of this class are $1, 2, 4, 9, 20, 45, 103, 235, 538, 1233, \ldots$.
	
	Let us recursively define the set $X_n$ of non-decreasing subexcedant functions over $[n]$, where $f\in X_n$ is obtained from $f'\in X_{n-1}$ by appending $n-1$ or $n$ at its end, or from $f\in X_{n-3}$ by appending the subword $(n-3)(n-3)(n-2)$ at its end. Thus, by construction we have the recurrence relation
	$$|X_n|=2|X_{n-1}|+|X_{n-3}|, n\geq 3, ~|X_0|=0, ~|X_1|=1, |X_2|=2.$$
	\begin{proposition}
		If $f\in X_n$, then $\sigma=\phi(f)\in\mathfrak{S}_n^{\scriptscriptstyle\nearrow}(231)$.
	\end{proposition}
	\begin{proof}
		Let $f\in X_n$ and $\sigma=\phi(f)$. Then $f$ is obtained from $f'\in X_{n-1}$ or $f'\in X_{n-3}$ by the operations defined above. Let $f\in X_n$ be obtained from $f'\in X_{n-1}$ by appending $n-1$ or $n$ at its end. Then $\sigma=(n,n-1)\circ\sigma'$ or $\sigma=(n,n)\circ\sigma'$, respectively, where $\sigma'=\phi(f')\in\mathfrak{S}_{n-1}^{\scriptscriptstyle\nearrow}(231)$. We show that $\sigma\in\mathfrak{S}_n^{\scriptscriptstyle\nearrow}(231)$. The case where $\sigma=(n,n)\circ\sigma'$ is obvious. Let $\sigma=(n,n-1)\circ\sigma'$ and assume that the pattern $231$ occurs in $\sigma$, i.\ e., there are integers $i<j<k$ such that $\sigma(j)>\sigma(i)>\sigma(k)$. If $\sigma(j)\neq n$, then since $\sigma(n)=n-1$ we have that $\sigma(i)=\sigma'(i), \sigma(j)=\sigma'(j), \sigma(k)=\sigma'(k)$, and thus the subword $\sigma'(i)\sigma'(j)\sigma'(k)$ would create the pattern $231$ in $\sigma'$. But this is a contradiction. If $\sigma(j)=n$, then since $\sigma'$ can be obtained from $\sigma$ by replacing the integer $n$ by $n-1$, $\sigma'(j)=n-1$ which is the largest integer in $\sigma'$. This implies that $\sigma'$ contains the pattern $231$, but this is also a contradiction. Therefore, $\sigma$ avoids the pattern $231$.
		
		Now consider that $f$ is obtained from $f'\in X_{n-3}$ by appending the subword $(n-3)(n-3)(n-2)$ at its end. Then we have $\sigma=(n,n-2)(n-1,n-3)(n-2,n-3)\sigma'$ and $\sigma(n-2)=n-1, \sigma(n-1)=n-3$ and $\sigma(n)=n-2$. If $\sigma$ contains the pattern $231$, then there exist integers $i<j<k\leq n-3$ with $\sigma(j)>\sigma(i)>\sigma(k)$. The operation of obtaining $\sigma'$ from $\sigma$ affects the pattern $231$ only if $\sigma(j)=n$. In this case we have $\sigma'(j)=n-3$ which is the largest integer in $\sigma'$. Thus, $\sigma'$ contains the pattern $231$  and this is a contradiction. 
	\end{proof}
	So $|\mathfrak{S}_n^{\scriptscriptstyle\nearrow}(231)|$ grows at least as fast as $|X_n|$.
	\subsection*{The case $312$.}
	In this case we show that $|\mathfrak{S}_n^{\scriptscriptstyle\nearrow}(312)|=2^{n-1}, n\ge1$.
	Let $Y_n:=\{f=f_1f_2\cdots f_n\in F_n^{\scriptscriptstyle\nearrow} : f_i=i, \forall i\in\Image(f)\}$. For instance, $Y_3=\{111, 113$, $122, 123\}$.
	\begin{proposition}
		A permutation $\sigma\in\mathfrak{S}_n^{\scriptscriptstyle\nearrow}(312)$  if and only if its transposition  array $f\in Y_n$.
	\end{proposition}
	\begin{proof}
		Assume that $f=f_1f_2\cdots f_n\in Y_n$ and $\sigma=\phi(f)$. Let $\Image(f)=\{\ell_1, \ell_2, \ldots, \ell_k\}$. Then $f_{\ell_i}=\ell_i, ~i\in \{1, 2, \ldots, k\}$. We know that $\ell_i$ is a fixed point of $f$ if and only if $\ell_i$ is the minimum of its cycle in $\sigma$. So all integers in the cycle of $\ell_i$ have $\nom$ equal to $\ell_i$,  hence they have to appear in increasing order in the cycle (when we write the cycle with its minimum $\ell_i$ at the beginning). Furthermore, since $f$ is a non-decreasing, all the integers in the cycle of $\ell_i$ form an integer interval, therefore  this cycle is $(\ell_i, ~\ell_i + 1, \ldots, \ell_{i+1}{-}1)$. 
		Thus, if $\sigma$ contains a $312$-pattern, i.\ e., there are $i<j<k$ such that $\sigma(i)>\sigma(k)>\sigma(j)$, then $i,j,k,\sigma(i), \sigma(j), \sigma(k)$ must be in the same cycle, say $C_r=(\ell_r,\sigma(\ell_r), \ldots)$, where $\ell_r<\sigma(\ell_r)<\sigma^2(\ell_r)<\cdots$. But this a contradiction. 
		
		Conversely, let $\sigma\in \mathfrak{S}_n^{\scriptscriptstyle\nearrow}(312)$ and $f=\phi^{-1}(\sigma)=f_1f_2\cdots f_n$. Suppose that $\Image(f)=\{\ell_1, \ell_2, \ldots, \ell_k\}$ and $\ell_j$ is the smallest image value of $f$ such that $f_{\ell_j}<\ell_j$. Let $t$ and $s$ be the positions of the leftmost and the rightmost occurrences of the value $\ell_j$ in $f$. Since $\ell_{j-1}\in\Image(f)$ is a fixed point it begins a cycle in $\sigma$, say $C$. Note that $\nom(\ell_j)=\ell_{j-1}$. Thus, $\ell_j, t, s\in C$. That is, $C=(\ell_{j-1}, \ell_{j-1}+1, \ldots, \ell_j{-}1, \ldots, t, \ldots, ~s, \ell_j, \ldots,~ t{-}1)$, where the integers between $\ell_j{-}1$ and $t$ (if any) are greater than $t$. Thus, there is a $312$ pattern in the positions $\ell_j{-}1$, $t{-}1$ and $s$, but this is a contradiction.
	\end{proof}
	\begin{corollary}
		The number of $312$-avoiding permutations in $\mathfrak{S}_n^{\scriptscriptstyle\nearrow}$ having exactly $k$ cycles is equal to $$|\{f\in Y_n : \IMA(f)=k\}|=\binom{n-1}{k-1},~ 0\leq k-1\leq n-1.$$
	\end{corollary}
	\begin{proof}
		For every subexcedant function $f=f_1f_2\cdots f_n$ we have $f_1=1$. Thus there are $n-1$ remaining positions from which we can choose $k-1$ positions, this can be possible in $\binom{n-1}{k-1}$, so that we have $f_i=i$. Then we fill in the remaining positions in one way to get the required function $f$. The number of cycles of the permutation coincides with $\IMA(f)$ because $f_i=i$ if and only if $i=\min(C)$, where $C$ is a cycle of $\sigma=\phi(f)$.
	\end{proof}
	\subsection*{The case $321$.} 
	In this case we can  provide a lower bound that shows that the size of this class grows exponentially with respect to $n$. The first few terms of this class are $1, 2, 5, 13, 35, 95, 261, 719, 1990, \ldots$.
	
	We note that a  permutation $\sigma\in\mathfrak{S}_n$ is $321$-avoiding if and only if both the subword of excedance letters and the  subword of anti-excedance letters  are increasing.  Therefore, a  permutation $\sigma\in\mathfrak{S}_n^{\scriptscriptstyle\nearrow}$ is $321$-avoiding if and only if the subword of excedance letters is increasing, i.\ e., if $i_1<i_2<\cdots <i_k$ is the increasing sequence of excedances of $\sigma$, then $\sigma(i_1)<\sigma(i_2)<\cdots<\sigma(i_k)$.
	Obviously, we obtain transposition arrays of  permutations satisfying this property  when we append $n$ or $n{-}1$ to the transposition array of a permutation  $\sigma$ in $\mathfrak{S}_{n-1}^{\scriptscriptstyle\nearrow}(321)$, therefore we have the following.
	\begin{proposition}
		For $n\ge1$, we have $$|\mathfrak{S}_{n}^{\scriptscriptstyle\nearrow}(321)|\geq 2^{n-1}.$$	
	\end{proposition} 
	\subsection*{Acknowledgement}
	The first author is grateful for the financial support extended by the cooperation agreement between International Science Program (ISP) at Uppsala University and Addis Ababa University. We appreciate the hospitality we got from IRIF during the research visits of the first author. We also thank J\"orgen Backelin of Stockholm University for the crucial discussions we had with him and his very useful suggestions.

\end{document}